\documentclass[twoside,leqno]{article}

\usepackage[letterpaper]{geometry}

\usepackage{siamproceedings}

\usepackage[T1]{fontenc}
\usepackage{amsfonts}
\usepackage{amssymb}
\usepackage{esint}
\usepackage{graphicx}
\usepackage{epstopdf}
\usepackage{enumitem}
\usepackage{algorithmic}
\ifpdf
  \DeclareGraphicsExtensions{.eps,.pdf,.png,.jpg}
\else
  \DeclareGraphicsExtensions{.eps}
\fi

\newsiamremark{remark}{Remark}
\newsiamremark{hypothesis}{Hypothesis}
\crefname{hypothesis}{Hypothesis}{Hypotheses}
\newsiamthm{claim}{Claim}
\newsiamthm{assumption}{Assumption}

\usepackage{amsopn}

\DeclareMathOperator{\tr}{tr}

\DeclareMathOperator{\dom}{dom}
\DeclareMathOperator{\ran}{ran}

\newcommand{\ud}{\,\mathrm{d}}

\newcommand{\RR}{\mathbb{R}}

\newcommand{\Or}{\mathcal{O}}

\newcommand{\wt}[1]{\widetilde{#1}}

\IfFileExists{mathabx.sty}%
  {\DeclareFontFamily{U}{mathx}{\hyphenchar\font45}%
   \DeclareFontShape{U}{mathx}{m}{n}{<->mathx10}{}%
   \DeclareSymbolFont{mathx}{U}{mathx}{m}{n}%
   \DeclareFontSubstitution{U}{mathx}{m}{n}%
   \DeclareMathAccent{\widebar}{0}{mathx}{"73}%
}{%
  \PackageWarning{mathabx}{%
    Package mathabx not available, therefore\MessageBreak substituting
    widebar with overline\MessageBreak }%
  \newcommand{\widebar}[1]{\overline{#1}}%
}

\newcommand{\mc}[1]{\mathcal{#1}}

\newcommand{\abs}[1]{\lvert#1\rvert}

\newcommand{\norm}[1]{\lVert#1\rVert}

\newcommand{\average}[1]{\langle#1\rangle}

\renewcommand{\Re}{\mathfrak{Re}}

\usepackage{todonotes}

\begin{document}

\newcommand\relatedversion{}

\title{\Large Quantitative Hypocoercivity and Lifting of Classical and Quantum Dynamics \relatedversion}
    
\author{Jianfeng Lu \thanks{Duke University  (\email{jianfeng.lu@duke.edu}).}}

\date{}

\maketitle


\fancyfoot[R]{\scriptsize{Copyright \textcopyright\ 20XX by SIAM\\
Unauthorized reproduction of this article is prohibited}}





\begin{abstract} We consider quantitative convergence analysis for hypocoercive dynamics such as Langevin and Lindblad equations describing classical and quantum open systems. Our goal is to provide an overview of recent results of hypocoercivity estimates based on space-time Poincar\'e inequality, providing a unified treatment for classical and quantum dynamics. Furthermore, we also present a unified lifting framework for accelerating both classical and quantum Markov semigroups, which leads to upper and lower bounds of convergence rates. 
\end{abstract}

\section{Introduction.}

Understanding and characterizing long-time behavior of dynamical systems have been an important field of study with important applications to various branches of mathematics, science, and engineering. For systems with dissipation, the dynamics is expected to converge to its equilibrium, regardless of whether the system is governed by classical or quantum equations of motion. 

In this work, we consider quantitative convergence analysis for a class of dynamics with degenerate dissipation, under the framework of hypocoercivity. Prototypical examples of such dynamics include Langevin dynamics (see \eqref{eq:langevin}) and Lindblad dynamics (see \eqref{eq:lindblad}) which characterize the evolution of classical and quantum open systems, respectively. The dissipation of such dynamics comes from interactions of the open systems with their environments, and it is common that such dissipation is degenerate: for example, the dissipation might act only on part of the system (such as the boundary of an open system) or only on part of the degrees of freedom (such as momentum variables). The hypocoercivity analysis aims to establish convergence to equilibrium by exploring the interplay between the conservative and dissipative parts of the dynamics. 

The study of hypocoercive dynamics has a long history, while we will focus on recent progress in quantitative hypocoercivity analysis based on a $L^2$ variational framework, initiated by \cite{albritton2024variational} and generalized and extended to various settings \cite{cao2023explicit, lu2022explicit, brigati2024explicit, brigati2025construct,li2024quantum,an2025convergence}. 
For both classical and quantum dynamics (Langevin and Lindblad), we obtain quantitative hypocoercive estimate with explicit convergence rate (\Cref{thm:classicalhypocoercive} and \Cref{thm:quantumhypocoercive}).
The key insight in this framework is to establish Poincar\'e type inequalities with augmented state space to take into account the effect of time evolution of the dynamics, such that the degenerate dissipation can propagate to the whole system. These inequalities are given as space-time Poincar\'e inequalities (\Cref{thm:spacetimepoincare} and \Cref{thm:quantumpoincare}) and flow Poincar\'e inequalities (\Cref{thm:flowpoincare}). Although all these results are available in recent papers by the author and collaborators \cite{cao2023explicit, li2024quantum, li2025speeding}, the hope is that this gives a more unified presentation of analysis for classical and quantum dynamics. 

As a consequence of the hypocoercivity estimates, with a suitable choice of the damping parameter, the dynamics exhibit accelerated convergence to equilibrium compared to their overdamped counterparts. The acceleration can be quadratic, known as the diffusive-to-ballistic dynamics transition in physics terms. This motivates the question to what extent hypocoercive dynamics can be used to accelerate convergence, which is analogous to the construction of lifted Markov chains to accelerate convergence in discrete state-space Markov chains \cite{diaconis2000analysis, chen1999lifting}. The recently developed framework of second-order lifting \cite{eberle2024non} proves to be useful to obtain upper and lower bounds of convergence rate for accelerated dynamics, which can be generalized to cover both classical and quantum dynamics \cite{li2025speeding}. In particular, we establish that the quadratic acceleration is optimal, as suggested by the generic upper bound of the convergence rate in terms of the spectral gap of the overdamped limit (\Cref{thm:upperbound}). Combined with the lower bound of convergence rate from hypocoercivity estimate, the lifting framework provides guidance to design optimal lifts for acceleration. 

\section{Langevin and Lindblad dynamics}

Let us first recall the Langevin dynamics for $(X_t, V_t) \in \RR^d \times \RR^d$
\begin{subequations}\label{eq:langevin}
  \begin{align}
    & \ud X_t = V_t \ud t; \\
    & \ud V_t = -\nabla U(X_t) \ud t - \gamma V_t \ud t + \sqrt{2\gamma} \ud W_t, 
  \end{align}
\end{subequations}
where $U: \RR^d \to \RR$ is a potential function that $U(x) \to \infty$ as $\abs{x} \to \infty$. The Langevin dynamics model particle systems with potential energy $U$ that interact with the environment: $X_t$ denotes the positions while $V_t$ denotes the corresponding momenta (for ease of notation, we have assumed that all particles have unit mass). The term $- \gamma V_t \ud t$ corresponds to damping, where $\gamma > 0$ is a friction parameter, and the term $\sqrt{2 \gamma} \ud W_t$ characterizes the random forcing from the environment given by the $d$-dimensional Wiener process $W_t$. The coefficients of these two terms are taken such that the system satisfies the fluctuation-dissipation relation. The Langevin equation is a prototypical model for classical open systems. 
We refer to the book \cite{pavliotis2014stochastic} for more background information on the Langevin equation. 

The law of Langevin dynamics, denoted as $\varrho(t, x, v)$, satisfies the kinetic Fokker-Planck equation
\begin{equation}\label{eq:kineticfokkerplanck}
    \partial_t \varrho = \bigl( - v \cdot \nabla_x + \nabla_x U \cdot \nabla_v \bigr) \varrho + \gamma \bigl( \nabla_v \cdot ( v \varrho) + \Delta_v \varrho \bigr). 
\end{equation}
Defining the Hamiltonian $h(x, v) = \frac{1}{2} \abs{v}^2 + U(x)$, we can verify that a stationary solution to \eqref{eq:kineticfokkerplanck} is given by
  \begin{equation}
    \mu(x, v) = \frac{1}{Z} e^{-h(x, v)}, \qquad \text{with normalizing constant } Z = \int e^{-h(x, v)} \ud x \ud v. 
  \end{equation}
  Under mild assumptions of $U$, $\mu$ is the density of the invariant measure of the Langevin dynamics \eqref{eq:langevin}. With some ambiguity of notation, we will also use the same notation for the invariant measure $\ud \mu = \mu(x, v) \ud x \ud v$. Note that $\mu$ is a product measure: $\ud \mu(x, v) = \ud \mu_x(x) \ud \kappa(v)$, where $\ud \mu_x \propto e^{-U(x)} \ud x$ and $\ud \kappa$ is the standard multivariate normal distribution. 
  
  We may also consider the backward Kolmogorov equation 
  \begin{equation}\label{eq:kolmogorov}
    \partial_t f = \mc{L} f := \bigl( v \cdot \nabla_x - \nabla_x U \cdot \nabla_v\bigr) f + \gamma \bigl(-v \cdot \nabla_v + \Delta_v\bigr) f. 
  \end{equation}
  This is equivalent to the kinetic Fokker-Planck equation \eqref{eq:kineticfokkerplanck}, as the latter can be written as 
  \begin{equation}
      \partial_t \varrho = \mc{L}^{\dagger} \varrho, 
  \end{equation}
  where $\mc{L}^{\dagger}$ denotes the $L^2$-adjoint of $\mc{L}$. Equivalent to the ergodicity of Langevin dynamics, the solution of the backward Kolmogorov equation as $t \to \infty$ converges to a constant as $\ker(\mc{L}) = \operatorname{span}\{1\}$. Note that $\int f(t, \cdot) \ud \mu$ is invariant under evolution \eqref{eq:kolmogorov}:
  \begin{equation}
    \partial_t \int f \ud \mu = \int \mc{L} f \ud \mu = \int f (\mc{L}^{\dagger} \mu) \ud x\ud v = 0.
  \end{equation} 
  Thus $f(t, x, v) \to \int f(0, \cdot) \ud \mu$ as $t \to \infty$.

We observe that the infinitesimal generator $\mc{L}$ consists of two components: 
    \begin{subequations}\label{eq:Ldecomp}
        \begin{align}
        & \mc{L} = \mc{L}_a + \gamma \mc{L}_s \qquad \text{with} \\
        & \mc{L}_a := v \cdot \nabla_x - \nabla_x U \cdot \nabla_v; \qquad \text{and} \\
        & \mc{L}_s := -v \cdot \nabla_v + \Delta_v.  
        \end{align}
    \end{subequations}
    Here $\mc{L}_a$ is the contribution from the Hamiltonian flow (with the Poisson bracket structure $\mc{L}_a = \frac{\partial h}{\partial v} \cdot\frac{\partial}{\partial x} - \frac{\partial h}{\partial x} \cdot \frac{\partial}{\partial v}$) and $\mc{L}_s$ takes into account fluctuation and dissipation from the interaction with the environment.
  
    Define the weighted $L^2$ inner product
    \begin{equation}
        \average{f, g}_{L^2(\mu)} = \int f(x, v) g(x, v) \ud \mu,  
    \end{equation}
    and the associated norm $\norm{\cdot}_{L^2(\mu)}$, we can check that $\mc{L}_a$ is antisymmetric and $\mc{L}_s$ is symmetric with respect to the weighted inner product: $\forall\, f, g$
    \begin{equation}\label{eq:symmetryantisymmetric}
        \average{f, \mc{L}_a g}_{L^2(\mu)} = - \average{\mc{L}_a f, g}_{L^2(\mu)} \qquad \text{and} \qquad \average{f, \mc{L}_s g}_{L^2(\mu)} = \average{\mc{L}_s f, g}_{L^2(\mu)}.
    \end{equation}  
    Note that the kernel of $\mc{L}_s$ is much larger than $\ker(\mc{L})$ as it contains all functions that are independent of $v$
    \begin{equation}\label{eq:kerLs}
        \ker(\mc{L}_s) = \{ f \mid \nabla_v f(x, v) = 0 \}.  
    \end{equation}

    \medskip
  
  For open quantum systems, under Markovian and semigroup assumptions, the evolution of the density matrix $\rho$ is given by the Lindblad equation \cite{lindblad1976generators, gorini1976completely}:
  \begin{equation}\label{eq:statelindblad}
    \begin{aligned}
      \partial_t \rho & = - i[H, \rho] + \gamma \sum_j \Bigl( \bigl[L_j \rho, L_j^{\dagger}\bigr] + \bigl[L_j, \rho L_j^{\dagger}\bigr] \Bigr) \\
      & = - i[H, \rho] + \gamma \sum_j \Bigl( L_j \rho L_j^{\dagger} - \frac{1}{2} \bigl\{ L_j^{\dagger} L_j, \rho \bigr\} \Bigr),
    \end{aligned}
  \end{equation}
  where $[\cdot, \cdot]$ and $\{ \cdot, \cdot\}$ are commutators and anti-commutators respectively: For two operators $A$, $B$ 
  \begin{equation*}
    [A, B] = AB - BA, \qquad \{A, B\} = AB + BA. 
  \end{equation*}
  Without the second term on the right-hand side of \eqref{eq:statelindblad}, we get the von Neumann equation $\partial_t \rho = - i [H, \rho]$ which characterizes the time evolution, in terms of the density matrix, of closed quantum system with the Hamiltonian operator $H$ (with a scaling choice so that the reduced Planck constant is $1$); the second term models the interactions of the system with environment, with $L_j$ being jump operators and $\gamma$ a parameter capturing the strength of the interaction. To avoid technicalities and keep the presentation simple, here we assume that the Hilbert space is finite dimensional and all involved operators are bounded. Generalizations will be considered in \S\ref{sec:Hilbertlifting}.
  
In the Heisenberg picture, considering the evolution of the observable, we have the following equivalent form of the Lindblad equation, analogous to the backward Kolmogorov equation.
  \begin{equation}\label{eq:lindblad}
    \partial_t A = \mc{L} A := i[H, A] + \gamma \sum_j \Bigl( L_j^{\dagger} [A, L_j] + [L_j^{\dagger}, A] L_j \Bigr).
  \end{equation}  
  We observe that the Lindblad equation has a rather similar structure as the backward Kolmogorov equation \eqref{eq:kolmogorov}: the first term of the right-hand sides of both equations are conservative, given by either the classical or quantum Hamiltonian as Poisson bracket or commutator, while the second term of the right-hand sides are dissipative modeling the interaction of the system with environment. As the classical Newtonian dynamics can be derived from Schr\"odinger equations in the semiclassical limit, the Langevin dynamics also arises in the semiclassical limit of Lindblad equations, see e.g., \cite{cao2017lindblad, hernandez2025classical, galkowski2025classical}. While we will not consider the semiclassical limit, the similarity of the structure of these equations is crucial, in fact one of the motivations of the study is to develop a unified theory understanding the convergence of classical and quantum dynamics for open systems.

  For a quantum state $\sigma$ (semi-positive definite operators with trace $1$) that is full rank (i.e., $\sigma$ as an operator is invertible), we define the Kubo-Martin-Schwinger (KMS) inner product  as 
  \begin{equation}\label{eq:KMS}
    \average{A, B}_{\sigma} := \tr( \sigma^{1/2} A^{\dagger} \sigma^{1/2} B) 
  \end{equation}
  and the associated norm 
  \begin{equation}
      \norm{A}_{\sigma} = \tr \bigl(\sigma^{1/2} A^{\dagger} \sigma^{1/2} A \bigr)^{1/2}.
  \end{equation}
  We say that the Lindbladian $\mc{L}$ satisfies the $\sigma$-KMS detailed balance condition if $\mc{L}$ is self-adjoint with respect to the inner product \eqref{eq:KMS} (for quantum dynamics, the notion of detailed balance is not unique as different inner product can be used, see e.g. \cite{temme2010chi, carlen2017gradient} for more discussions).
  \begin{equation}
      \average{A, \mc{L} B}_{\sigma} = \average{\mc{L} A, B}_{\sigma}, \quad \forall\, A, B.  
  \end{equation}
  In particular, taking $B$ to be the identity operator, we get $\tr(A, \mc{L}^{\dagger} \sigma) = 0$, for any $A$, and thus $\mc{L}^{\dagger} \sigma = 0$ so that $\sigma$ is a stationary state of the Lindblad evolution \eqref{eq:statelindblad}. 

  Similarly to \eqref{eq:Ldecomp}, we observe that the generator consists of two parts: $\mc{L} = \mc{L}_a + \gamma \mc{L}_s$ with
  \begin{subequations}\label{eq:Qdecomp}
    \begin{align}
      & \mc{L}_a = i [H, \cdot] \\
      & \mc{L}_s = \sum_j \Bigl( L_j^{\dagger} [\cdot, L_j] + [L_j^{\dagger}, \cdot] L_j \Bigr).
    \end{align}
  \end{subequations}
  We will assume the following symmetry and anti-symmetry hold for $\mc{L}_s$ and $\mc{L}_a$. \begin{equation}\label{eq:qunantumsymmetryantisymmetric}
        \average{A, \mc{L}_a B}_{\sigma} = - \average{\mc{L}_a A, B}_{\sigma} \qquad \text{and} \qquad \average{A, \mc{L}_s B}_{\sigma} = \average{B, \mc{L}_s A}_{\sigma}.
    \end{equation}  
    It follows that $\sigma$ is a stationary state of $\mc{L}$. 

Moreover, we observe that the symmetric part of the generator may have a non-trivial kernel (cf.~\eqref{eq:kerLs}): 
\begin{equation}
    \ker(\mc{L}_s) = \{ A \mid [L_j, A] = 0,\, [L_j^{\dagger}, A] = 0 \}.
\end{equation}
In addition, we have the following characterization for $\ker(\mc{L})$ (see \cite[Theorem 7.2]{wolf2012quantum})
\begin{equation}
    \ker(\mc{L}) =  \{ A \mid [H, A] = 0,\, [L_j, A] = 0,\, [L_j^{\dagger}, A] = 0 \}. 
\end{equation}
Thus in general, $\ker(\mc{L}) \subsetneq \ker(\mc{L}_s)$: $\ker(\mc{L})$ is a strict subspace of $\ker(\mc{L}_s)$, which again is analogous to the case of Langevin dynamics. 

\section{Hypocoercivity and space-time Poincar\'e inequality} 
\label{sec:spacetimepoincare}

Our primary focus is the convergence of Langevin and Lindblad dynamics to equilibrium. 
To motivate the discussion, let us consider a cousin of the Langevin dynamics, the overdamped dynamics given by 
\begin{equation}\label{eq:overdamp}
    \ud X_t = - \nabla U(X_t) \ud t + \sqrt{2} \ud W_t. 
\end{equation}
In fact, as will be discussed in \Cref{sec:overdamp}, it is the asymptotic limit of the Langevin dynamics when $\gamma \to \infty$. Analogous to \eqref{eq:kolmogorov}, we consider the associated backward Kolmogorov equation 
\begin{equation}\label{eq:overdampkolmogorov}
    \partial_t f(t, x) = (- \nabla_x U \cdot \nabla_x + \Delta_x) f(t, x) =: \mc{L}_O f(t, x) 
\end{equation}
which converges to a constant as $t \to \infty$; equivalently, $\mu_x \propto e^{-U(x)}$ is the invariant measure of the overdamped dynamics \eqref{eq:overdamp}. Note that $\int f(t, \cdot) \ud \mu_x$ is invariant under evolution \eqref{eq:overdampkolmogorov}. Taking an initial condition $f(0, x)$ such that $\int f(0, x) \ud \mu_x = 0$, to quantify the convergence, we consider the weighted $L^2$ norm and get 
\begin{equation}\label{eq:energy}
    \begin{aligned}
        \frac{\ud}{\ud t} \int \abs{f(t, x)}^2 \ud \mu_x 
        & = 2 \int f(t, x) \mc{L}_O f(t, x) \ud \mu_x \\
        & = - 2 \int \abs{\nabla_x f(t, x)}^2 \ud \mu_x,
    \end{aligned}
\end{equation}
where the last equality is a result of integration by parts. If the right hand side can be bounded above by $- \norm{f(t, x)}_{L^2(\mu_x)}^2$, we arrive at exponential decay of the $L^2(\mu_x)$-norm of the solution. This motivates the following. 
\begin{assumption}\label{assump:poincare}
    The measure $\mu_x = \frac{1}{Z_x} e^{-U(x)}$ satisfies the Poincar\'e inequality: For any $f(x)$ such that $\int f \ud \mu_x = 0$, we have 
    \begin{equation}\label{eq:poincare}
        m \int \abs{f(x)}^2 \ud \mu_x \leq \int \abs{\nabla f(x)}^2 \ud \mu_x = \average{f, (-\mc{L}_O) f}_{L^2(\mu_x)}
    \end{equation}
\end{assumption}
Note that the Poincar\'e inequality can be viewed as the coercivity of the generator $-\mc{L}_O$ with respect to the $L^2(\mu_x)$ norm. Combining  \eqref{eq:energy} with Assumption~\ref{assump:poincare}, we have
\begin{equation}
    \frac{\ud}{\ud t} \int \abs{f(t, x)}^2 \ud \mu_x \leq -2m \int \abs{f(t, x)}^2 \ud \mu_x,
\end{equation}
and thus 
\begin{equation}\label{eq:overdampexpdecay}
    \norm{f(t, \cdot)}_{L^2(\mu_x)} \leq e^{-mt} \norm{f(0, \cdot)}_{L^2(\mu_x)}.
\end{equation}

It is natural to ask if a similar argument can be carried out for the underdamped Langevin dynamics, mimicking \eqref{eq:poincare}, we consider the Dirichlet form for $f = f(x, p)$
\begin{equation}
    \average{f, (- \mc{L}) f}_{L^2(\mu)} = \average{f, (- \mc{L}_a) f}_{L^2(\mu)} + \gamma \average{f, (-\mc{L}_s) f}_{L^2(\mu)}. 
\end{equation}
The first term on the right hand side vanishes as $\mc{L}_a$ is antisymmetric (recall \eqref{eq:symmetryantisymmetric}), and an integration by parts in $v$ gives 
\begin{equation}
    \average{f, (- \mc{L}_s) f}_{L^2(\mu)} = \int \abs{\nabla_v f(x, v)}^2 \ud \mu.
\end{equation}
Therefore, $-\mc{L}_s$ cannot be coercive as it is impossible that $\norm{\nabla_v f(x, v)}_{L^2(\mu)}$ can control the $L^2(\mu)$ norm of $f$. This is of course due to the fact that $\mc{L}_s$ has a non-trivial kernel as seen in \eqref{eq:kerLs}. 

\medskip 

Nevertheless, for the (underdamped) Langevin dynamics \eqref{eq:langevin}, we can still establish exponential convergence to equilibrium with quantitative rate estimates, as shown in the following theorem. Such results are referred as quantitative hypocoercive estimates. The study of hypocoercivity in the context of Langevin type equations has a long history dated back to Kolmogorov \cite{kolmogoroff1934zufallige} and H\"ormander \cite{hormander1967hypoelliptic}. Quantitative estimates have been established more recently by Villani~\cite{villani2009hypocoercivity} using a $H^1$-hypocoercivity framework. More explicit estimates of the convergence rate are achieved by a $L^2$-approach developed by Dolbeault, Mouhot and Schmeiser \cite{dolbeault2009hypocoercivity, dolbeault2015hypocoercivity} (see also \cite{herau2006hypocoercivity}). Our estimate is based on a variational framework initiated by \cite{albritton2024variational} and further developed and extended to various settings \cite{cao2023explicit, lu2022explicit, brigati2024explicit, brigati2025construct,li2024quantum,an2025convergence}. This approach establishes the following hypocoercive estimate for Langevin dynamics.

\begin{theorem}[Hypocoercivity of Langevin dynamics \cite{cao2023explicit}] \label{thm:classicalhypocoercive}
Assume $U$ is convex, superlinear as $\abs{x} \to \infty$: $\lim_{x\to \infty} \frac{U(x)}{\abs{x}^{\beta}} \to \infty$ for some $\beta > 1$, and satisfy $\abs{\nabla^2 U} \lesssim 1 + \abs{\nabla U}$. 
Let $f(t, x, v)$ solve the backward Kolmogorov equation \eqref{eq:kolmogorov} corresponding to the Langevin dynamics, we have 
\begin{equation}
    \biggl\lVert f(t, \cdot) - \int f(t, \cdot) \ud \mu \biggr\rVert_{L^2(\mu)} \leq C e^{-\nu t} \,\biggl\lVert f(0, \cdot) - \int f(0, \cdot) \ud \mu \biggr\rVert_{L^2(\mu)}.
\end{equation}
The convergence rate $\nu$ is given explicitly by 
\begin{equation}
    \nu = \frac{m \gamma}{c (\sqrt{m} + \gamma)^2 },
\end{equation}
where $c$ and $C$ are some universal constants. 
\end{theorem}

Taking the friction coefficient $\gamma$ to be $\sqrt{m}$, we have $\nu = \sqrt{m} / (4c)$, and thus for $f$ satisfying $\int f \ud \mu = 0$, we have 
\begin{equation}\label{eq:expdecay}
    \norm{f(t, \cdot)}_{L^2(\mu)} \leq C e^{-\sqrt{m}t / (4c)} \norm{f(0, \cdot)}_{L^2(\mu)}.
\end{equation}
Compared to \eqref{eq:overdampexpdecay}, \eqref{eq:expdecay} shows that Langevin dynamics exhibit quadratic acceleration of convergence when $m \ll 1$. This is similar to convex optimization, where underdamped type algorithms, such as conjugate gradient, heavy ball, and Nesterov's algorithm achieve quadratic acceleration compared to steepest descent in various settings \cite{nesterov1983method}. In fact, as the overdamped dynamics is a Wasserstein gradient flow in the probability space \cite{jordan1998variational}, one may try to push the analogy further, see e.g., \cite{ma2021there, chewi2025logconcave}. Furthermore, we can also verify explicitly in the case of $U$ being a quadratic function that the $\Theta(\sqrt{m})$-convergence rate estimate is sharp. The sharpness can be also understood from the lifting point of view \cite{eberle2024non}, as will be elaborated in \Cref{sec:lifting}. From a physics point of view, such acceleration corresponds to a diffusive-to-ballistic transition of the dynamics. We also remark that the exponential decay of hypocoercive dynamics as \eqref{eq:expdecay} would necessarily have a prefactor $C$ that is strictly larger than $1$; in fact, if exponential decay holds with some rate $\nu$ and $C = 1$, it can be shown that the generator has to be coercive.

Let us make some comments on the assumptions of \Cref{thm:classicalhypocoercive}: 
\begin{itemize}
\item The convexity assumption of $U$ can be slightly relaxed, see \cite{cao2023explicit} for more general cases; however, it is not expected that quadratic acceleration can be achieved when the measure $\mu_x$ exhibits metastability, i.e., the potential $U(x)$ contains several local minima separated by large potential barriers. The convergence rate in that case is given by the Eyring-Kramers formula, see e.g., \cite{bovier2004metastability, nier2005hypoelliptic}. 
\item The super-linearity assumption of $U$ can be relaxed to allow distributions with heavier tails, see \cite{brigati2024explicit}. 
\end{itemize}

The key to obtain \Cref{thm:classicalhypocoercive} is a space-time Poincar\'e inequality, which plays a similar role in energy estimate as the standard Poincar\'e inequality (Assumption~\ref{assump:poincare}). However, as the generator $\mc{L}$ lacks coercivity, the dissipation effects in the $x$-variable can only be obtained taking time evolution into consideration, thus the space-time Poincare inequality considers a time average over an interval $[0, T]$.  

\begin{theorem}[Space-time Poincar\'e inequality \cite{cao2023explicit}] \label{thm:spacetimepoincare}
    Under the same assumptions of \Cref{thm:classicalhypocoercive}, we have 
    \begin{multline}\label{eq:spacetimepoincare}
        \Biggl( \fint_0^T \Bigl \lVert f - \fint_0^T \int f \ud \mu \ud t   \Bigl\rVert_{L^2(\mu)}^2 \ud t \Biggr)^{1/2} \lesssim 
        \Bigl(1 + \frac{1}{T \sqrt{m}} \Bigr) \Biggl( \fint_0^T \norm{(I - \Pi_v)  f}_{L^2(\mu)}^2 \ud t \Biggr)^{1/2} \\ 
        + \Bigl(\frac{1}{\sqrt{m}} + T \Bigr) \Biggl( \fint_0^T \norm{ ( I - \mc{L}_s)^{-1/2} (\partial_t - \mc{L}_a) f }_{L^2(\mu)}^2 \ud t \Biggr)^{1/2},
    \end{multline}
    where $\fint_0^T = \frac{1}{T} \int_0^T$ is the averaged time integral. Here $\Pi_v$ is the orthogonal projection to $\ker(\mc{L}_s)$ with respect to the inner product $\average{\cdot, \cdot}_{L^2(\mu)}$, given by 
    \begin{equation}
        (\Pi_v f)(t, x) := \int f(t, x, v) \ud \kappa(v).
    \end{equation}
\end{theorem}

\begin{remark}
Since the $v$-marginal of $\mu$ is given by standard Gaussian $\kappa$, by Gaussian Poincar\'e inequality, we have 
\begin{equation}\label{eq:GaussianPoincare}
    \norm{(I - \Pi_v) f}_{L^2(\mu)}^2 \leq \norm{\nabla_v f}_{L^2(\mu)}^2 = \average{f, (-\mc{L}_s) f}_{L^2(\mu)}
\end{equation}
and thus the first term on the right-hand side of \eqref{eq:spacetimepoincare} can be replaced by 
\begin{equation*}
\Bigl(1 + \frac{1}{T \sqrt{m}} \Bigr) \Biggl( \fint_0^T \norm{\nabla_v  f}_{L^2(\mu)}^2 \ud t \Biggr)^{1/2},
\end{equation*}
which resembles more the usual Poincar\'e inequality. In addition, if $f$ is a solution to \eqref{eq:kolmogorov}, we have $(\partial_t - \mc{L}_a)f = \gamma \mc{L}_s f$, and thus the second term on the right-hand side of \eqref{eq:spacetimepoincare} is also equivalent to the dissipation term under the underdamped dynamics. 
\end{remark}

\medskip 

As Lindblad dynamics can be viewed as a quantum generalization of Langevin dynamics, it is thus natural to ask whether hypocoercivity estimates can be extended to the evolution of open quantum systems. Such estimate for quantum dynamics was first established in \cite{fang2025mixing} using the $L^2$-hypocoercivity approach similar to \cite{dolbeault2009hypocoercivity, dolbeault2015hypocoercivity} for classical linear kinetic equations. 
The estimate based on a quantum extension of the space-time Poincar\'e inequality (\Cref{thm:quantumpoincare}) was subsequently established in \cite{li2024quantum}.

For hypocoercivity estimates, we require some structural assumptions of the Lindblad dynamics \eqref{eq:lindblad}. Recall that we have already assumed the symmetry and anti-symmetry of $\mc{L}_s$ and $\mc{L}_a$ in \eqref{eq:qunantumsymmetryantisymmetric} and thus it follows that $\sigma$ is a stationary state. 
Denote by $\lambda(\mc{L}_s)$ the spectral gap of $\mc{L}_s$, which we assume to be strictly positive: 
\begin{equation}
    \lambda(\mc{L}_s) := \inf_{A \in \ker(\mc{L}_s)^{\perp}} \frac{\norm{\mc{L}_s A}_{\sigma}}{\norm{A}_{\sigma}} > 0. 
\end{equation}
Since $\mc{L}_s$ is self-adjoint with respect to the KMS inner product $\average{\cdot, \cdot}_{\sigma}$, we have, analogous to \eqref{eq:GaussianPoincare}, 
\begin{equation}
    \lambda(\mc{L}_s) \norm{ (I - \Pi_s) A}_{\sigma}^2 \leq \average{A, (- \mc{L}_s) A}_{\sigma},
\end{equation}
where $\Pi_s$ denotes the orthogonal projection to $\ker(\mc{L}_s)$ with respect to the KMS inner product $\average{\cdot, \cdot}_{\sigma}$.

\begin{theorem}[Quantum space-time Poincar\'e inequality \cite{li2024quantum}]  \label{thm:quantumpoincare}
Assume the above structural assumptions and that $\Pi_s \mc{L}_a \Pi_s = 0$, we have 
    \begin{multline}
    \Biggl( \fint_0^T \Bigl \lVert A(t) - \fint_0^T \tr(\sigma A(t)) \ud t  \Bigl\rVert_{\sigma}^2 \ud t \Biggr)^{1/2} \leq
            C_{1, T} \Biggl( \fint_0^T \norm{(I - \Pi_s) A(t)}_{\sigma}^2   \ud t \Biggr)^{1/2} \\  
        + C_{2, T} \Biggl( \fint_0^T \norm{ (I - \mc{L}_s)^{-1/2} (\partial_t - \mc{L}_a ) A(t) }_{\sigma}^2 \ud t \Biggr)^{1/2},
    \end{multline}
    where the constants $C_{1, T}$ and $C_{2, T}$ can be made explicit (see \cite{li2024quantum} for details).
\end{theorem}

The role of the assumption $\Pi_s \mc{L}_a \Pi_s = 0$ will be further elaborated in \Cref{sec:lifting}, here let us just comment that an analogous condition is automatically satisfied in the case of Langevin dynamics. Using the quantum space-time Poincar\'e inequality in an energy estimate, we obtain the hypocoercive estimates for Lindblad dynamics. 
\begin{theorem}[Hypocoercivity of Lindblad dynamics \cite{li2024quantum}]\label{thm:quantumhypocoercive}
Under the same assumptions as \Cref{thm:quantumpoincare}, let $A(t)$ solve the Lindblad equation \eqref{eq:lindblad}, we have 
\begin{equation}\label{eq:Adecay}
    \norm{A(t) - \tr(\sigma A(t))}_{\sigma} \leq C  e^{-\nu t} \norm{A(0) - \tr(\sigma A(0))}_{\sigma},
\end{equation}
with prefactor $C = e^{\nu T}$ and convergence rate
\begin{equation}
    \nu = \frac{\gamma \lambda(\mc{L}_s)}{C_{1, T}^2 + \gamma^2 \lambda(\mc{L}_s) C_{2, T}^2},
\end{equation}
where $C_{1, T}$ and $C_{2, T}$ are the same constants in \Cref{thm:quantumpoincare}. Here, the parameter $T$ can be optimized to achieve larger convergence rate $\nu$ in \eqref{eq:Adecay}. 
\end{theorem}

\section{Lifting}\label{sec:lifting} 

Compared with the overdamped dynamics, the results discussed in the previous section indicate that by extending the dynamics to include the additional momentum variable, we might achieve accelerated convergence. The idea of constructing faster converging Markov processes in a larger state space has been previously studied for discrete Markov chains \cite{diaconis2000analysis, chen1999lifting, vucelja2016lifting} under the notion of lifting. From this perspective, we can view the Langevin dynamics as a lifting of the overdamped dynamics. This concept is formalized by \cite{eberle2024non} and is tightly connected with hypocoercivity (see e.g., \cite{brigati2024hypocoercivity, eberle2025convergence, li2025speeding}), which will be discussed in this section. We will first recall the overdamped limit in \S\ref{sec:overdamp}, which to some extent is reciprocal to lifting. In order to cover both classical and quantum dynamics, we will consider semigroups in Hilbert space, recalled in \S\ref{sec:Hilbertlifting}. The lifting for these semigroups will be defined in \S\ref{sec:deflifting}. The main advantage of the lifting point of view is that it provides a framework to establish lower and upper bounds of convergence rates of the lifted dynamics combined with hypocoercive estimates, which we will discuss in \S\ref{sec:convergence}.

\subsection{Overdamped limit}\label{sec:overdamp}

Recall that in Langevin dynamics \eqref{eq:langevin}, $\gamma$ is the friction parameter; the quantitative hypocoercivity estimate \cref{thm:classicalhypocoercive} suggests the optimal choice $\gamma = \Theta(\sqrt{m})$ for convergence. 
Instead, if we consider the limit $\gamma \to \infty$ with a rescaling of the time $t \mapsto t / \gamma$, we arrive at the overdamped Langevin dynamics (see e.g., \cite[Chapter 6.5]{pavliotis2014stochastic})
\begin{equation}
    \partial_t \varrho(t, x) = \nabla_x \cdot \bigl(\varrho \nabla_x U + \nabla_x \varrho \bigr)(t, x). 
\end{equation}
\Cref{thm:classicalhypocoercive} indicates that the Langevin dynamics can achieve a quadratic acceleration of convergence compared to its overdamped limit. 

In general, motivated by the structure of kinetic Fokker-Planck and Lindblad equations, we may consider a family of generators (cf. \eqref{eq:Ldecomp} and \eqref{eq:Qdecomp}). Here and in the sequel, we use the subscript $\gamma$ in $\mc{L}_{\gamma}$ to emphasize its dependence on the friction parameter. While for the corresponding semigroups, with some abuse of notation, we will still use $\mc{P}_t$ to keep the notation simple. 
\begin{equation}
    \mc{L}_{\gamma} = \mc{L}_a + \gamma \mc{L}_s, \quad \gamma > 0,
\end{equation}
where $\mc{L}_a$ and $\mc{L}_s$ are symmetric and anti-symmetric with respect to an appropriate inner product: 
\begin{equation}\label{eq:symmetry}
    \mc{L}_a^{\star} = - \mc{L}_a, \qquad \text{and} \qquad \mc{L}_s^{\star} = \mc{L}_s.
\end{equation}
Consider the evolution 
\begin{equation}
    \partial_t A^{\gamma}(t) = \bigl(\mc{L}_a + \gamma \mc{L}_s \bigr) A^{\gamma}(t)
\end{equation}
and the limit $\gamma \to \infty$. Assume that the solution admits a formal asymptotic expansion: 
\begin{equation*}
    A^{\gamma}(t) = A_0(t) + \gamma^{-1} A_1(t) + \gamma^{-2} A_2(t) +  \Or(\gamma^{-2}).
\end{equation*}
Substituting this expansion into the evolution equation yields
\begin{align*}
    \gamma \mc{L}_s A_0 + 
    \bigl( \mc{L}_s A_1 + \mc{L}_a A_0 - \partial_t A_0 \bigr) 
    + \Or(\gamma^{-1}) = 0.    
\end{align*}
Matching powers of $\gamma$ leads to:
\begin{subequations}
\begin{align}
     \Or(\gamma): \quad &\mc{L}_s A_0 = 0, \label{subeq1} \\ 
     \Or(1): \quad &\mc{L}_s A_1 = \partial_t A_0 - \mc{L}_a A_0. \label{subeq2} 
\end{align}
\end{subequations}
From \eqref{subeq1}, we deduce that $\Pi_s A_0(t) = A_0(t)$ for all $t \geq 0$, where $\Pi_s$ denotes the orthogonal projection on $\ker \mc{L}_s$. Combining this with \eqref{subeq2} leads to the effective dynamics:
\begin{align} \label{eq:eff1st}
    \partial_t A_0(t) = \Pi_s \mc{L}_a \Pi_s A_0(t).
\end{align}
Therefore, in summary, $A^{\gamma}(t)$ is approximated by $A_0(t) \in \ker(\mc{L}_s)$, governed by \eqref{eq:eff1st}. 

However, when $\Pi_s \mc{L}_a \Pi_s = 0$, the right-hand side of the leading-order dynamics \eqref{eq:eff1st} vanishes. To obtain non-trivial dynamics, we need to look at longer time horizon. Thus, let us rescale time $t \mapsto t/\gamma$ and consider
\begin{align} \label{eq:asym2}
    \partial_t A^{\gamma}(t) = (\gamma \mc{L}_a + \gamma^2 \mc{L}_s) A^{\gamma}(t).
\end{align}
Proceeding again using the formal asymptotic expansion and matching orders, we have 
\begin{subequations}
\begin{align}
     \Or(\gamma^2): \quad &\mc{L}_s A_0 = 0, \label{eq:0order} \\
     \Or(\gamma^1): \quad &\mc{L}_s A_1  = -\mc{L}_a A_0, \label{eq:1order} \\
     \Or(1): \quad & \mc{L}_s A_2 = \partial_t A_0 - \mc{L}_a A_1. \label{eq:norder}
\end{align}
\end{subequations}
The leading order equation \eqref{eq:0order} implies $\Pi_s A_0 = A_0$ as before, while the solvability condition for \eqref{eq:1order} requires $\ran(\mc{L}_a \Pi_s) \subset \ran(\mc{L}_s)$. This is equivalent to the following structural assumption: 
\begin{assumption} \label{assump:PHP}
    $\Pi_s \mc{L}_a \Pi_s = 0$.
\end{assumption}
Recall that in the case of Langevin, $\mc{L}_a = v \cdot \nabla_x - \nabla_x U \cdot \nabla_v$ and $\Pi_s f = \int f(x, v) \ud \kappa(v)$, Assumption~\ref{assump:PHP} can be explicitly verified. In the case of Lindbladian, this assumption is a key structural assumption in \Cref{thm:quantumpoincare}. 
Under Assumption~\ref{assump:PHP}, we obtain the first-order correction:
\begin{equation*}
    A_1 \in -\mc{L}_s^{-1} \mc{L}_a A_0 + \ker(\mc{L}_s),
\end{equation*}
where $\mc{L}_s^{-1}$ denotes the pseudoinverse of $\mc{L}_s$. To close the equation for the leading-order term $A_0$, we apply the projection $\Pi_s$ to both sides of \eqref{eq:norder}, which gives:
\begin{equation} \label{eq:eff2nd}
    \begin{aligned}
        \partial_t A_0 &= \Pi_s \mc{L}_a A_1 \\
        &= -\Pi_s\mc{L}_a \mc{L}_s^{-1} \mc{L}_a A_0 \\
        &= - (\mc{L}_a \Pi_s)^\star (-\mc{L}_s)^{-1} (\mc{L}_a \Pi_s) A_0,
    \end{aligned}
\end{equation}
where in the last equality we have used the symmetry of $\mc{L}_a$ and $\mc{L}_s$ \eqref{eq:symmetry}. Thus, in this regime, the dynamics $A^{\gamma}(t)$ is effectively described by $A_0(t) \in \ker(\mc{L}_s)$ evolved as \eqref{eq:eff2nd} with an effective generator given by 
\begin{equation} \label{eq:odlimit_generator}
    \mc{L}_O := - (\mc{L}_a \Pi_s)^{\star} (- \mc{L}_s)^{-1} (\mc{L}_a \Pi_s).
\end{equation}
In the case of Langevin dynamics, we can explicitly check that the obtained $\mc{L}_O$ coincides with the infinitesimal generator of the overdamped dynamics. Similar limiting dynamics has also been considered for open quantum systems; see e.g., \cite{zanardi2014coherent, li2024quantum}. 

In the general case, it is natural to consider the convergence rate of semigroup generated by $\mc{L}_{\gamma}$ and its connection with the overdamped limit governed by $\mc{L}_O$. Alternatively, given a generator $\mc{L}_O$ we may construct generator $\mc{L}_{\gamma}$ with optimal choice of $\gamma$ in the hope of accelerating convergence. This motivates the framework of lifting, proposed for classical Markov process in \cite{eberle2024non}. In the following, we will discuss a generalization of lifting to semigroups on Hilbert spaces, and thus covering both classical and quantum dynamics.

\subsection{Semigroups on Hilbert spaces}\label{sec:Hilbertlifting}

Let $\{\mc{P}_t\}_{t \ge 0}$ be a contractive strongly continuous semigroup on a Hilbert space $\mc{H}$, in other words, $\mc{P}_t$ satisfies the properties: (i) $\mc{P}_0 = \mathrm{Id}$; (ii) $\mc{P}_s \mc{P}_t = \mc{P}_{s + t}$ for $s, t \ge 0$; (iii) $\lim_{t \to 0^+}\norm{\mc{P}_t x - x}_{\mc{H}} = 0$ for all $x \in \mc{H}$; (iv) $\norm{\mc{P}_t}_{\mc{H} \to \mc{H}} \leq 1$. Conditions (i)--(iii) define a $C_0$-semigroup and condition (iv) imposes contractivity. 

The generator $\mc{L}$ corresponding to the semigroup $\{\mc{P}_t\}_{t \geq 0}$ is a closed, densely defined operator given by
\begin{align*}
    \mc{L} x := \lim_{t \to 0^+} \frac{\mc{P}_t x - x}{t},
\end{align*}
with $\dom(\mc{L})$ consisting of all $x \in \mc{H}$ for which the limit exists in the norm topology of $\mc{H}$. In particular, unlike the finite dimensional case, here and in the sequel, we allow the generator $\mc{L}$ to be unbounded. The following is standard, following from the Hille-Yosida and Lumer-Phillips theorems \cite{engel2000one}.

\begin{lemma} \label{lem:dissipative}
Let $\mc{P}_t$ be a $C_0$-semigroup with generator $(\mc{L},\,\dom(\mc{L}))$. It holds that:
\begin{itemize}
    \item $x_t := \mc{P}_t x_0$ for $x_0 \in \dom(\mc{L})$ is continuous differentiable on $[0, \infty)$ and satisfies $\dot{x}_t = \mc{L} x_t$ with initial condition $x_{t = 0} = x_0$, and $x_t \in \dom(\mc{L})$ for all $t \ge 0$, that is, $\mc{P}_t (\dom (\mc{L})) \subset \dom(\mc{L}).$
    
    \item If $\mc{P}_t$ is contractive,  we have  $ \sigma(\mc{L}) \subset \{\lambda \in \mathbb{C} \mid \Re \lambda \le 0\}$ and $\norm{(\lambda - \mc{L})^{-1}}_{\mc{H}} \le (\Re \lambda)^{-1}$ for $\lambda \in \mathbb{C}$ with $\Re \lambda > 0$. Moreover, $\mc{L}$ is dissipative: 
    \begin{align*}
    \Re \average{ x, \mc{L} x}_{\mc{H}} \le 0\,,\quad \forall\, x \in \dom(\mc{L})\,.
    \end{align*}
\end{itemize}
\end{lemma}

The equilibrium subspace of the semigroup $\{\mc{P}_t\}_{t \geq 0}$ is given by the kernel of $\mc{L}$
\begin{align} \label{eq:mixsubspace}
    \{x \in \mc{H} \mid \mc{P}_t x = x \text{ for all } t \geq 0\} = \ker(\mc{L}).
\end{align}
We denote by $\mc{P}_{\infty}$ the orthogonal projection on $\ker{\mc{L}}$. 
We say that $\{\mc{P}_t\}_{t \geq 0}$ \emph{converges to equilibrium} if for each $x \in \mc{H}$, there exists a decay function $r: \mathbb{R}_+ \to \mathbb{R}_+$ with $\lim_{t \to \infty} r(t) = 0$ such that
\begin{align} \label{eq:exponentialconverge}
    \|\mc{P}_t x - \mc{P}_{\infty} x\|_{\mc{H}} \leq r(t) \|x - \mc{P}_{\infty} x\|_{\mc{H}}.
\end{align}
When $r(t) = Ce^{- \nu t}$ for some $C >  1, \nu > 0$, we say the semigroup is \emph{hypocoercive}. In the case of $r(t) = e^{-\nu t}$ (i.e., $C = 1$), the semigroup is \emph{coercive}. 
    For either case, the sharp convergence rate 
\begin{equation}\label{eq:defspectralgap}
    \nu_0 := \sup \bigl\{ \nu > 0 \mid \exists\, C \geq 1 \text{ such that }    \norm{\mc{P}_t x - \mc{P}_{\infty} x}_{\mc{H}} \leq C e^{-\nu t} \norm{x - \mc{P}_{\infty} x}_{\mc{H}}, \forall\, x\in \mc{H}  \bigr\}
\end{equation}
is given by the spectral gap of $\mc{L}$ \cite[Chapter IV]{engel2000one}
\begin{equation}
    \nu_0 = \lambda(\mc{L}) := \inf \bigl\{\Re(\lambda) \mid \lambda \in \operatorname{spec}(-\mc{L}) \backslash \{0\} \bigr\}. 
\end{equation}
Although this gives the asymptotic rate as $t \to \infty$, often in applications, we care about the time such that $\mc{P}_t$ almost reaches equilibrium, which motivates the definition of the relaxation time as 
\begin{equation}
    t_{\text{rel}}(\mc{L}) := \inf 
    \bigl\{ t \geq 0 \mid \norm{\mc{P}_t x - \mc{P}_{\infty}x }_{\mc{H}} \leq e^{-1} \norm{x - \mc{P}_{\infty} x}_{\mc{H}}, \forall\, x \in \mc{H} \bigr\}.
\end{equation}
As shown in \cite{eberle2024non, li2024quantum}, the relaxation time can be lower bounded using the singular value gap of $\mc{L}$, given by the spectral gap of $\abs{\mc{L}} = \sqrt{\mc{L}^{\star}\mc{L}}$, 
\begin{equation}\label{eq:trelbound}
    t_{\text{rel}}(\mc{L}) \geq \frac{1}{2\, \mathfrak{s}(\mc{L})}.
\end{equation}
Note that in general $\mc{L}$ is not self-adjoint, and thus its singular value $\mathfrak{s}(\mc{L})$ is not necessarily the spectral gap as in \eqref{eq:defspectralgap}. 
We have a perhaps slightly more explicit formula for the singular value gap: 
\begin{equation}
    \mathfrak{s}(\mc{L}) = \inf\{\norm{\mc{L} x}_{\mc{H}} \mid x \in \dom(\mc{L}) \cap \ker{\mc{L}}^\perp,\, \norm{x}_{\mc{H}} = 1\}. \label{eq:singulargap}
\end{equation}

From the definition of the relaxation time and hypocoerivice semigroups, we have the following lemma. 
\begin{lemma} \label{lem:singulargapgeneral}
    Let $\mc{P}_t$ be a hypocoercive $C_0$-semigroup with generator $\mc{L}$, with decay function $r(t) = C e^{-\nu t}$ in \eqref{eq:exponentialconverge}, it holds 
    \begin{align*}
        \nu \leq (1 + \log C)\, \mathfrak{s}(\mc{L}).
    \end{align*}
\end{lemma}
We remark that the singular value gap is also used to analyze relaxation speed of non-reversible discrete Markov chains \cite{chatterjee2025spectral}. 

\subsection{Lifting in Hilbert spaces}\label{sec:deflifting}

Let us now formalize lifting in Hilbert spaces, which can be understood as the opposite of the overdamped limit. 
A hypocoercive semigroup $\mc{P}_t$ can be viewed as a \emph{lifting} of a symmetric coercive semigroup $\mc{P}_{t,O}$ acting on a closed strict subspace, denoted as $\mc{H}_O$, with the induced inner product $\average{x, y}_{\mc{H}_O} = \average{x,y}_{\mc{H}}$ for $x,y \in \mc{H}_O$. For the semigroup $\mc{P}_{t, O}$, referred as the \emph{collapsed semigroup}, we assume that it is coercive on $\mc{H}_{O}$ with rate $\lambda_O > 0$:
\begin{align*}
    \|\mc{P}_{t,O} x\|_{\mc{H}_{O}} \leq e^{-\lambda_O t} \|x \|_{\mc{H}_{O}}, 
\end{align*}
where $\lambda_O$ is the spectral gap of $\mc{L}_O$. 
We now give the formal definition of lifting, which extends the concept of the second-order lift of classical Markov processes \cite{eberle2024non,eberle2025convergence}.

\begin{definition}[Lifting] \label{def:generallift}
Let $\mc{P}_t$ and $\mc{P}_{t,O}$ be contraction $C_0$-semigroups on $\mc{H}$ and $\mc{H}_{O}$ with generators $\mc{L}$ and $\mc{L}_{O}$, respectively. Assume $\mc{H}_O \subset \mc{H}$ and $\dom(\mc{L}_{O}) \subset \dom(\mc{L})$. 
 $\mc{P}_t$ is a \emph{second-order lifted semigroup} of $\mc{P}_{t,O}$ if
 \begin{subequations}     
    \begin{enumerate}[label=(\roman*)]
        \item \label{condg2}  There holds, for any $x \in \mc{H}_{O}$, $y \in \dom(\mc{L}) \cap \mc{H}_{O}$, 
        \begin{align} \label{eq:phpgen}
            \average{x, \mc{L} y}_{\mc{H}} = 0 \,. 
        \end{align}
        \item \label{condg3} There exists a positive  bounded operator $\mc{S}: \mc{H}_{O}^\perp \to  \mc{H}_{O}^\perp$ such that
        \begin{align} \label{eq:liftgen}
            \average{\mc{L} x,  \mc{S} \mc{L} y}_{\mc{H}} = - \average{ x, \mc{L}_{O} y}_{\mc{H}_{O}}, \quad  \forall x \in \dom(\mc{L}) \cap \mc{H}_{O},\, y \in \dom(\mc{L}_{O}). 
        \end{align}
    Note that \eqref{eq:phpgen} implies $\mc{L} y \in \mc{H}_{O}^\perp$, and therefore $\mc{S}\mc{L} y$ is well defined. 
    \end{enumerate}
    \end{subequations}
\end{definition}
We include an operator $\mc{S}$ to make the framework more flexible. Recall that in the overdamped limit, $\mc{S}$ is given by $(-\mc{L}_s)^{-1}$, see \eqref{eq:odlimit_generator}. The definition of lifting used in \cite{eberle2024non,eberle2025convergence} for classical Markov process corresponds to the choice $\mc{S} = I$. 

While we focus on second-order lifting, it is also possible to consider first-order lifting, so that the generator $\mc{L}$ satisfies 
\begin{align}\label{eq:first_order_condition}
    \langle x, \mc{L} y \rangle_{\mc{H}} = - \langle x, \mc{L}_{O} y \rangle_{\mc{H}_{O}}, \quad \forall x,y \in \dom(\mc{L}_{O})\,.
\end{align}
instead of conditions \eqref{eq:phpgen} and \eqref{eq:liftgen} in \Cref{def:generallift}. The first order lifting is useful  in the context of accelerating discrete-time finite-state Markov chains \cite{chen1999lifting}; see e.g., discussions in \cite[Remark 12]{eberle2024non}.

\subsection{Convergence rate of lifted dynamics}\label{sec:convergence}
An immediate implication of lifting in \cref{def:generallift} and \cref{lem:singulargapgeneral} is that the convergence rate of the semigroup $\mc{P}_t$ generated by $\mc{L}$ is at most $\mathcal{O}(\sqrt{\lambda_O})$ (\cref{thm:upperbound}), where $\lambda_O$ is the spectral gap corresponding to the semigroup $\mc{P}_{t, O}$; i.e., lifting can lead to at most a quadratic speed-up. 

\begin{theorem}[Upper bound of convergence rate \cite{li2025speeding}] \label{thm:upperbound}
    Suppose $\mc{P}_t$ is a second-order lift of $\mc{P}_{t, O}$ and the constant $\wt{s}_{\rm m} := \inf_{x \in \mc{H}_O^\perp \backslash \{0\}} \frac{\norm{\Pi_1 \mc{S} \Pi_1 x}_{\mc{H}_O}}{\norm{x}_{\mc{H}_O}}$ is positive, where $\Pi_1$ is the projection from $\mc{H}_O^\perp$ to $\overline{\ran(\mc{L}|_{\mc{H}_O})}$. 
    Then, there holds
    \begin{equation} \label{eq:upperrate}
        \nu \leq (1 + \log C) \sqrt{(\wt{s}_{\rm m})^{-1}\lambda_O}\,.
    \end{equation}
\end{theorem}

The intuition behind this upper bound is that $\mc{L}_O$ can be viewed as roughly the square of the lifted generator $\mc{L}^{\star} \mc{L}$ due to the relation \eqref{eq:liftgen}; thus the square of the singular value gap of $\mc{L}$ cannot be larger than the spectral gap of $\mc{L}_O$. 

\smallskip 

The hope would be then to construct lifting that achieves the quadratic acceleration, which requires estimates on the lower bound of the convergence rate. Unlike the upper bound, such estimates would require more structures of the generator; in what follows, we will limit ourselves to a family of generators $\mc{L}_{\gamma}$ with structures similar to the Langevin and Lindblad dynamics. While the upper bound in \Cref{thm:upperbound} is independent of the parameter $\gamma$, the lower bound would suggest choice of $\gamma$ to achieve optimal convergence.

\begin{assumption} \label{assmp:decomposition}  
Assume that $(\mc{L}_s,\dom(\mc{L}_s))$ and $(\mc{L}_a,\dom(\mc{L}_a))$ are closed and densely defined operators on $\mc{H}$. Consider the family of generators 
\begin{align*}
    \mc{L}_{\gamma}  = \mc{L}_a + \gamma \mc{L}_s, \quad \gamma > 0,
\end{align*}
with $\dom(\mc{L}) = \dom(\mc{L}_s)\cap \dom(\mc{L}_a)$. We assume
\begin{enumerate} [label=(\roman*)]
    \begin{subequations}
        \item \label{condsym} $\mc{L}_s$ is symmetric and satisfies 
        \begin{equation} \label{eq:kernells}
            \Pi \mc{L}_s x = 0\,, \quad \forall\, x \in \dom(\mc{L}), 
        \end{equation}
        where $\Pi$ is the orthogonal projection from $\mc{H}$ to $\mc{H}_O$. 
        \item $\mc{L}_s$ is coercive: For some $\lambda_S > 0$,
        \begin{align} \label{eq:lscoercive}
        \average{x, (-\mc{L}_s) x}_{\mc{H}} \geq 
            \lambda_S \norm{x - \Pi x}_{\mc{H}}^2\,, \quad\forall\, x \in \dom(\mc{L}).  
        \end{align} 
        \item \label{condlift}  $\mc{L}_{a}$ is a lift of $\mc{L}_{O}$, in particular, $\dom(\mc{L}_O) \subset \dom(\mc{L}_a)$.   
        \item \label{condantisym} $\mc{L}_a$ is anti-symmetric on $\dom(\mc{L}_{O})$:
        \begin{align} \label{eq:weakanti}
           \mc{L}_a^* x = - \mc{L}_a x\,, \quad\forall\, x \in \dom(\mc{L}_{O}). 
        \end{align}
        \item\label{condker} There holds 
        \begin{align} \label{eq:kernelconst}
            \ker(\mc{L}_{O}) = \ker(\mc{L}_{\gamma}) \subsetneq \mc{H}_O.     
        \end{align}
        Thus the equilibrium space of $\mc{P}_t$ and $\mc{P}_{t, O}$ is the same strict subspace of $\mc{H}_O$. 
    \end{subequations}
\end{enumerate}
\end{assumption}
Let us remark that under the domain assumption $\mc{H}_O \subset \dom(\mc{L}_s)$, it can verified that $\mc{L}_a$ is a lift of $\mc{L}_{O}$ if and only if $\mc{L}_{\gamma}$ is a lift of $\mc{L}_{O}$ for any $\gamma > 0$, which explains Assumption~\ref{assmp:decomposition}\ref{condlift}. In particular, this applies to finite-dimensional dynamics, since \eqref{eq:kernells} and \eqref{eq:lscoercive} imply $\ker(\mc{L}_s) = \mc{H}_O$.

The following lemma characterizes the equilibria of $\mc{P}_t$ in terms of $\ker(\mc{L}_a)$ and $\ker(\mc{L}_s)$. It  justifies Assumption~\ref{assmp:decomposition}\ref{condker}  which implicitly indicates that equilibrium subspace $\ker(\mc{L}_{\gamma})$ of $\mc{P}_t$ is independent of $\gamma$. 

\begin{lemma} \label{lem:kernelrela}
    Let $\mc{P}_t$ be a hypocoercive $C_0$-semigroups on $\mc{H}$ satisfying Assumption~\ref{assmp:decomposition} conditions \ref{condsym}--\ref{condantisym}. Then, $x \in \dom (\mc{L})$ is an equilibrium, i.e., $\mc{L}_{\gamma} x = 0$, if and only if $\mc{L}_s x = 0$ and $\mc{L}_ax = 0$.
\end{lemma}

\begin{remark}\label{rem:hypocoerpt}
If instead of being a strict subspace as in Assumption~\ref{assmp:decomposition}\ref{condker}, $\ker(\mc{L}_{\gamma}) = \mc{H}_O$, then $\mc{P}_t$ is in fact coercive. This is because for $x \in \dom(\mc{L})$, we have 
\begin{align} \label{eq:hypocoercive11}
    \Re \average{x, (-\mc{L}_{\gamma}) x}_{\mc{H}} = \gamma \average{x, (-\mc{L}_s) x}_{\mc{H}} \geq \lambda_S \gamma \norm{x - \Pi x}_{\mc{H}}^2\,.
\end{align} 
If $\ker{\mc{L}_{\gamma}} = \mc{H}_O$, the projection $\Pi$ to the space $\mc{H}_O$ then coincides with $\mc{P}_{\infty}$, the projection to $\ker(\mc{L}_{\gamma})$, thus the above inequality \eqref{eq:hypocoercive11} gives the coercivity of $\mc{L}_{\gamma}$: $\Re \average{x, (-\mc{L}_{\gamma}) x}_{\mc{H}} \geq \lambda_S \gamma \norm{x - \mc{P}_\infty x}^2_{\mc{H}}$. 
\end{remark}

\medskip 

Our goal is to establish quantitative hypocoercivity for $\mc{P}_t$
\begin{equation}
    \norm{\mc{P}_t x_0}_{\mc{H}} \leq C e^{-\nu t} \norm{x_0}_{\mc{H}}.
\end{equation}  
As before, the general strategy is to control $\norm{x}_{\mc{H}}$, as in classical Poincar\'e type inequalities. 
To proceed, we start by decomposing the norm using the orthogonal projection $\Pi$ from $\mc{H}$ to $\mc{H}_O$: 
\begin{align} \label{eq:flowdecomposition}
    \norm{x}_{\mc{H}}^2 = \norm{\Pi x}_{\mc{H}}^2 + \norm{x - \Pi x}_{\mc{H}}^2.
\end{align}
The term $\norm{x - \Pi x}_{\mc{H}}^2$ can be immediately bounded by the coercivity assumption of $\mc{L}_s$ \eqref{eq:lscoercive}.
Controlling $\norm{\Pi x}_{\mc{H}}^2$ on the other hand is not trivial: For $\Pi x \in \dom(\mc{L}_{O}) \cap \dom(\mc{L})$, we have 
\begin{align*}
    \Re \langle \Pi x, (-\mc{L}) \Pi x \rangle_{\mc{H}} &= \gamma \langle \Pi x, (-\mc{L}_s) \Pi x \rangle_{\mc{H}} + \Re \langle \Pi x, (-\mc{L}_a)\Pi x \rangle_{\mc{H}} = 0\,,
\end{align*}
by \eqref{eq:kernells} and \eqref{eq:weakanti}, and thus the dissipation of $\mc{P}_t$ vanishes on $\mc{H}_O$. This calculation also explains why the condition \ref{condantisym} of Assumption~\ref{assmp:decomposition} imposes anti-symmetry of $\mc{L}_a$ only in the dense subspace $\dom(\mc{L}_{O}) \subset \mc{H}_O$. 

To address the lack of coercivity, the idea is to consider the evolution over an interval of time, so that the interplay between $\mc{L}_a$ and $\mc{L}_s$ would lead to effective dissipation, and the time-average behavior would exhibit coercivity. This already appears in our discussion in \Cref{sec:spacetimepoincare}, while in the current setup, we will use a slightly simplified framework based on the flow Poincaré inequality inspired by \cite{eberle2025convergence}. 
It is a simplified version of the space-time Poincaré inequalities, as instead of considering any possible $x_t \in L^2([0,T];\mc{H})$, the flow Poincar\'e inequality applies to the solution $x_t = \mc{P}_t x_0$.
It takes the form of standard Poincar\'e inequality (see Assumption~\ref{assump:poincare}), but augmented with a time variable: there exists $\alpha_T > 0$ depending on $T > 0$ such that for any $x_0 \in \dom(\mc{L}_{\gamma}) \cap \ker(\mc{L}_{\gamma})^\perp$,
\begin{align}\label{eq:flow_poincare}
    \alpha_T \fint_0^T \|x_t\|_{\mc{H}}^2 \ud t \leq \fint_0^T \langle x_t, (-\mc{L}_s) x_t \rangle_{\mc{H}} \ud t\,,  \quad \text{where } x_t := \mc{P}_t x_0.
\end{align}
More precisely, we have the following theorem for the hypocoercive dynamics $\mc{P}_t$. 

\begin{theorem}[Flow Poincar\'e inequality \cite{li2025speeding}] \label{thm:flowpoincare}
Let $\mc{P}_t$ be a hypocoercive $C_0$-semigroup on $\mc{H}$ with the equilibrium subspace $\ker(\mc{L}_{\gamma})$. Under Assumptions~\ref{assmp:decomposition} and some additional technical assumptions, for any time horizon $T > 0$, and any $x_0 \in \ker(\mc{L}_{\gamma})^\perp \cap \dom(\mc{L}_{\gamma})$, there holds: 
\begin{align} \label{eq:tspoincare}
  \frac{1}{T} \int_0^T \norm{x_t}_{\mc{H}}^2 \ud t \leq \bigl(C_{1, T} + \gamma^2 C_{2, T}\bigr) \frac{1}{T} \int_0^T \average{x_t, (-\mc{L}_s) x_t}_{\mc{H}} \ud t,
\end{align}
where $x_t: = \mc{P}_t x_0$ and $C_{1, T}$ and $C_{2, T}$ are some explicit constants independent of $\gamma$ (for details, see \cite{li2025speeding}).
\end{theorem}
The proof relies on the use of the lifting structure and technical a priori estimates for solutions to an abstract divergence equation, we will not go into the details and refer the interested readers to \cite{eberle2025convergence, li2025speeding}.
With the flow Poincar\'e inequality, a standard energy estimate gives exponential decay in the $L^2([0,T];\mc{H})$ norm, as follows.
\begin{theorem}[Lower bound of convergence rate \cite{li2025speeding}] \label{thm:lowerboundpt}
Under the same assumptions as in \cref{thm:flowpoincare}, it holds that for any period $T > 0$ and any initial $x_0 \in \ker(\mc{L})^\perp$, 
\begin{align} \label{eq:expconverge}
    \fint_t^{t + T} \norm{\mc{P}_s x_0}_{\mc{H}}^2 \ud s \le e^{- 2\nu t} \norm{x_0}_{\mc{H}}^2\,,
\end{align}
with explicit convergence rate 
\begin{equation} \label{eq:ratel2}
    \nu = \frac{\gamma}{C_{1, T} + \gamma^2 C_{2, T}},
\end{equation}
where $C_{1,T}$ and $C_{2,T}$ are the flow Poincar\'e constants in \cref{thm:flowpoincare}.
\end{theorem}

As a corollary, we obtain that for any $T > 0$ and $x_0 \in \ker(\mc{L})^\perp$, we have the hypocoercive estimate
\begin{equation}\label{eq:expcoverct}
    \norm{\mc{P}_t x_0}_{\mc{H}} \le e^{\nu T} e^{-\nu t} \norm{x_0}_{\mc{H}}.
\end{equation}

The convergence rate given in \eqref{eq:ratel2} is maximized by taking $\gamma_{\max} = \sqrt{C_{1, T}/C_{2, T}}$; since by the upper bound the convergence rate cannot exceed $\Or(\sqrt{\lambda_O})$, if the lower bound \eqref{eq:ratel2} from \Cref{thm:lowerboundpt} matches in order, this gives us the optimal lifted dynamics in terms of accelerated convergence. We refer the interested readers to \cite{li2025speeding} for examples of optimal lifted dynamics, which we will not go into details here.

\section{Discussions}

We have focused our discussion on hypocoercivity analysis for Langevin and Lindblad dynamics, similar approaches based on space-time Poincar\'e inequality can be applied to other dynamics with degenerate dissipation. In particular, the sampling dynamics of various piecewise deterministic Markov processes, such as bouncy particle methods \cite{bouchard2018bouncy} and zigzag sampler \cite{bierkens2019zig}, share similar structures as Langevin dynamics and can be analyzed using similar approaches
\cite{lu2022explicit}. 
In addition to sampling applications, hypocoercivity analysis can also be applied to dynamics that arises from min-max optimizations
\cite{wang2023local, an2025convergence}.

While our discussion has focused on continuous-time dynamics, time discretization is required in order to use Langevin dynamics as a Markov chain Monte Carlo sampling algorithm for the invariant measure. The quantitative convergence analysis of discretization of Langevin dynamics has received much attention in recent years and remains quite active, we refer to recent work \cite{altschuler2024faster, altschuler2025shifted} and references therein. 

For applications of dissipative quantum dynamics, several Lindblad dynamics have been proposed in recent literature for preparing thermal and ground states \cite{chen2023quantum, chen2023efficient, ding2024single, ding2025efficient}, mostly based on Lindblad dynamics with detailed balance. It is of interest to consider possible acceleration of such dynamics based on lifting.

\section*{Acknowledgments.}
This research has been supported in part by the National Science Foundation through the awards DMS-2012286 and DMS-2309378.  
We are grateful to our collaborators Jing An, Yu Cao, Di Fang, Bowen Li, Lihan Wang, and Yu Tong. We also thank Andreas Eberle, Jonathan Mattingly and Gabriel Stoltz for helpful discussions. 

\bibliographystyle{siamplain}
\bibliography{hypocoercivity}

\end{document}